\magnification=\magstep1
\parindent=0pt
\baselineskip=13pt

\font\BF=cmbx12 at 17pt
\font\bfsl=cmbxsl9
\font\slsl=cmsl9
\font\am=msam10
\font\amam=msam7
\font\bm=msbm10

\def\eop{\hfill{\am\char"03}}
\def\loe{\mathrel{\hbox{\am\char"36}}}
\def\lloe{\mathrel{\hbox{\amam\char"36}}}
\def\goe{\mathrel{\hbox{\am\char"3E}}}
\def\emptyset{\hbox{\bm\char"3F}}
\def\NN{\hbox{\bm\char"4E}}
\def\RN{\hbox{\bm\char"52}}
\setbox8=\hbox{\rm{\bf1)} \ }
\def\indent{${}\kern\wd8$}
\def\hb#1{\hbox to\wd8{{\bf#1)} \ }}
\setbox9=\hbox{\rm(O) \ }
\def\CA{$\hbox to\wd9{\rm(\hfill A\hfil) \ }$}
\def\CB{$\hbox to\wd9{\rm(\hfil B\hfil) \ }$}
\def\CH{$\hbox to\wd9{\rm(\hfil H\hfil) \ }$}
\def\CL{$\hbox to\wd9{\rm(\hfil\hfil L\hfil) \ }$}
\def\CO{$\hbox to\wd9{\rm(O) \ }$}
\def\Cr{$\hbox to\wd9{\rm(r)\hfil}$}
\def\abc{[\,a,b\,]\cap[\,s-\delta,s+\delta\,]}


\centerline{\BF Elementary real analysis}
\medskip
\centerline{\BF without compactness argument}

\bigskip\medskip


\hrule
\medskip
\centerline{Claude-Alain Faure}
\medskip
\hrule

\bigskip\medskip


{\baselineskip=11.7pt
{\bfsl Abstract.\ \ \slsl In a famous paper, R.\ A.\ Gordon proved a dozen
theorems using tagged partitions and Cousin's theorem. The purpose of this
paper is to present several classical results using the key-lemma underlying
Cousin's theorem.}
\par}

\bigskip


{\bf First key-lemma.} \ Let $\cal I$ be a set of subintervals of $[\,a,b\,]$. We
suppose that this set $\cal I$ satisfies the following conditions:
\smallskip
(A) \ $[\,a,s\,]\in\cal I$ and $[\,s,t\,]\in{\cal I}\ \Rightarrow\ [\,a,t\,]\in\cal I$\quad%
(additivity)
\smallskip
\CL for every $s\in{]\,a,b\,]}$, there exists $a\loe r<s$ such that $[\,x,s\,]\in\cal I$
for all $r\loe x<s$
\smallskip
\Cr for every $s\in{[\,a,b\,[}$, there exists $s<t\loe b$ with $[\,s,t\,]\in\cal I$
\smallskip
Then the whole interval $[\,a,b\,]$ must be an element of $\cal I$.

\medskip


{\bf Proof.} \ We consider the set $E=\{a\}\cup\{x\in{]\,a,b\,]}\mathrel/[\,a,x\,]\in
{\cal I} \,\}$ and $s=\sup(E)$.
\smallskip
We first show that $s\in E$. We may assume $s>a$. Let $a\loe r<s$ be as in
condition (L). By definition of the supremum, there exists $x\in{]\,r,s\,]}\cap E$.
If $x=s$, we are done. And if $x<s$, then $[\,a,x\,]\in\cal I$ and $[\,x,s\,]\in{\cal
I}$ $\Rightarrow$ $[\,a,s\,]\in\cal I$ by additivity. So $s\in E$.
\smallskip
We now show that $s=b$, which proves the lemma. Suppose on the contrary
that $s<b$. By condition (r) we have $[\,s,t\,]\in\cal I$ for some $s<t\loe b$,
and hence $[\,a,t\,]\in \cal I$ by additivity. So $t\in E$, in contradiction to the
supremacy of $s$.\eop

\bigskip


{\bf Remark.} This lemma was used in [5] to prove extreme value theorem
and intermediate value theorem. Its full strength applies in Theorem 4. But 
in a few special cases, such as Heine's theorem, the following result seems
more efficient.

\bigskip


{\bf Second key-lemma.} \ Let $\cal I$ be a set of subintervals of $[\,a,b\,]$.
We suppose that this set\break $\cal I$ satisfies the following conditions:
\smallskip
\CO$[\,a,s\,]\in\cal I$, $[\,r,t\,]\in{\cal I}$ and $r<s\loe t$ $\Rightarrow$ $[\,a,t\,]
\in\cal I$\quad(overlapping additivity)
\smallskip
\CB for every $s\in[\,a,b\,]$, there exists $\delta>0$ with $[\,r,t\,]={[\,a,b\,]}\cap
{[\,s-\delta,s+\delta\,]}\in\cal I$
\smallskip
Then the whole interval $[\,a,b\,]$ must be an element of $\cal I$.

\medskip


{\bf Proof.} \ We consider the set $E=\{x\in{]\,a,b\,]}\mathrel/[\,a,x\,]\in{\cal I}\,\}$.
By condition (B) with $s=a$, we conclude that $E$ is not empty. Therefore $s=
\sup(E)$ exists and $s>a$. Now let $[\,r,t\,]$ be as in condition (B). By definition
of the supremum, there exists $x\in{]\,r,s\,]}\cap E$. Then $[\,a,x\,]\in\cal I$, $[\,r,
t\,]\in\cal I$ and $r<x\loe s\loe t$ $\Rightarrow$ $[\,a,t\,]\in\cal I$ by additivity. So
$t\in E$, and hence $t\loe s$. We therefore obtain $s=t=b\in E$, which proves
the lemma.\eop

\vfill\eject


{\bf Theorem 1.} \ If $f:[\,a,b\,]\to\RN$ is locally bounded, then it is bounded on
$[\,a,b\,]$.

\medskip


{\bf Proof.} \ We consider the set ${\cal I}=\bigl\{I\subset{[\,a,b\,]}\mathrel{\bigm/}f$
is bounded on $I\bigr\}$.
\smallskip
\CO If $f$ is bounded on $[\,a,s\,]$ and $[\,r,t\,]$, then it is bounded on $[\,a,s\,]
\cup[\,r,t\,]$.
\smallskip
\CB By hypothesis, there exists $\delta>0$ such that $f$ is bounded on $\abc$.
\smallskip
By the second key-lemma, we conclude that $[\,a,b\,]\in\cal I$.\eop

\bigskip


{\bf Theorem 2.} \ Any zero-free continuous function $f:[\,a,b\,]\to\RN$ has constant
sign.

\medskip


{\bf Proof.} \ We consider the set ${\cal I}=\{I\subset[\,a,b\,]\mathrel/f$ has constant
sign on $I\}$.
\smallskip
\CO Let $[\,a,s\,]\in\cal I$ and $[\,r,t\,]\in\cal I$ with $r<s\loe t$. Since $[\,a,s\,]\cap
[\,r,t\,]\not=\emptyset$, the sign must be the same on both intervals. Then $[\,a,t\,]
\in\cal I$.
\smallskip
\CB Let $s\in[\,a,b\,]$. Say $f(s)>0$. By continuity, there exists $\delta>0$ such that
$f(x)>0$ for all $x\in\abc$. Then $\abc\in\cal I$.
\smallskip
By the second key-lemma we conclude that $[\,a,b\,]\in\cal I$.\eop

\bigskip


{\bf Theorem 3.} \ Any upper semicontinuous function $f:C\to\RN$ defined on a
(non-empty) closed subset  $C\subset[\,a,b\,]$ has a maximum.

\medskip


{\bf Proof.} \ Suppose on the contrary that $f$ has no maximum. We consider
the set
\medskip
\centerline{${\cal I}=\{I\subset[\,a,b\,]\mathrel/\hbox{there exists }y\in C$ such
that $f(x)<f(y)$ for all $x\in I\cap C\}$}
\medskip
\CO If $f(x)<f(y_1)$ for all $x\in[\,a,s\,]\cap C$ and if $f(x)<f(y_2)$ for all $x\in
[\,r,t\,]\cap C$, then $f(x)<\max{\{f(y_1),f(y_2)\}}$ for all $x\in\bigl([\,a,s\,]\cup
[\,r,t\,]\bigr)\cap C$.
\smallskip
\CB Let $s\in[\,a,b\,]$. {\sl Case 1:} $s\notin C$. There exists $\delta>0$ with
$[\,s-\delta,s+\delta\,]\cap C=\emptyset$. Then $\abc\in\cal I$. {\sl Case 2:}
$s\in C$. According to the assumption on $f$, there exists $y\in C$ with $f
(s)<f(y)$. And by semicontinuity, there exists $\delta>0$ such that $f(x)<f(y)$
for all $x\in[\,s-\delta,s+\delta\,]\cap C$. Then $\abc\in\cal I$.
\smallskip
By the second key-lemma we conclude that $[\,a,b\,]\in\cal I$, which is clearly
impossible.\eop

\bigskip


{\bf Theorem 4.} \ Let $f:[\,c,d\,]\to\RN$ satisfy the following conditions:
\smallskip
{\bf1)} \ for every $s\in{]\,c,d\,]}$, there exists $c\loe r<s$ such that $f(x)<f(s)$ for
all $r\loe x<s$
\smallskip
{\bf2)} \ for all $s\in{[\,c,d\,[}$ and $s<t\loe d$, there exists $s<x\loe t$ with $f(s)<
f(x)$
\smallskip
Then $f$ is strictly increasing on the interval $[\,c,d\,]$.

\medskip


{\bf Proof.} \ Let $c\loe a<b\loe d$. We consider the set ${\cal I}=\bigl\{[\,u,v\,]
\subset[\,a,b\,]\mathrel{\bigm/}f(u)<f(v)\bigr\}$.
\smallskip
\CA Let $[\,a,s\,]\in\cal I$ and $[\,s,t\,]\in\cal I$ . Then $f(a)<f(s)<f(t)$ $\Rightarrow$
$|\,a,t\,]\in\cal I$.
\smallskip
\CL Let $s\in{]\,a,b\,]}$. Then condition 1) implies that $[\,x,s\,]\in\cal I$ for all $\max
\{a,r\}\loe x<s$.
\smallskip
\Cr Let $s\in{[\,a,b\,[}$. Then condition 2) implies that there exists $s<x\loe b$ with
$[\,s,x\,]\in\cal I$.
\smallskip
By the first key-lemma we conclude that $[\,a,b\,]\in\cal I$, which proves that $f(a)<
f(b)$.\eop

\bigskip
\settabs=2\columns


{\bf Corollary 1.} \ Let $f:[\,c,d\,]\to\RN$ satisfy the following conditions:
\smallskip
\+ {\bf1)} \ $\underline{D\!}\,{}_-f(s)>0$ for all $s\in{]\,c,d\,]}$ & {\bf2)} \ $\,\overline
{\!D}{}_+f(s)>0$ for all $s\in{[\,c,d\,[}$ \cr
\smallskip
Then $f$ is strictly increasing on the interval $[\,c,d\,]$.

\medskip


{\bf Proof.} \ {\bf1)} \ If $\underline{D\!}\,{}_-f(s)=\sup\limits_{r<s}\,\inf\limits_{\vrule
width 0pt height 6.2pt r\lloe x<s}\displaystyle{f(s)-f(x)\over s-x}>0$, then condition
1) is satisfied.
\smallskip
{\bf2)} \ If $\,\overline{\!D}{}_+f(s)=\inf\limits_{\vrule width 0pt height 6.8pt s<t}\,\sup
\limits_{s<x\lloe t}\displaystyle{f(x)-f(s)\over x-s}>0$, then condition 2) is satisfied.
\eop

\bigskip


{\bf Corollary 2.} \ Let $f:[\,c,d\,]\to\RN$ satisfy the following conditions:
\smallskip
\+ {\bf1)} \ $\underline{D\!}\,{}_-f(s)\goe0$ for all $s\in{]\,c,d\,]}$ & {\bf2)} \ $\,
\overline{\!D}{}_+f(s)\goe0$ for all $s\in{[\,c,d\,[}$ \cr
\smallskip
Then $f$ is increasing on the interval $[\,c,d\,]$.

\medskip


\line{{\bf Proof.} \ Suppose on the contrary that $a<b$ and $f(a)>f(b)$. By a classic
argument, we}
\smallskip
consider the function $g(x)=f(x)-\,\displaystyle{f(b)-f(a)\over b-a}\,(x-a)$ to derive a
contradiction.\eop

\bigskip


{\bf Heine's theorem.} \ Let $f:C\to\RN$ be a continuous function defined on a
closed subset $C\subset[\,a,b\,]$. Then $f$ is uniformly continuous.

\medskip


{\bf Proof.} \ Let $\varepsilon>0$. We consider the set $\cal I$ of all subintervals
$I\subset[\,a,b\,]$ for which there exists $\delta>0$ such that $x,y\in I\cap C$ and
$|x-y|\loe\delta$ $\Rightarrow$ $|f(x)-f(y)|<\varepsilon$.
\smallskip
\CO Let $[\,a,s\,]\in\cal I$ and $[\,r,t\,]\in{\cal I}$ with $r<s\loe t$. There exist a
suitable $\delta{}_1>0$ for $[\,a,s\,]$ and a suitable $\delta{}_2>0$ for $[\,r,t\,]$.
Then $\delta=\min{\{\delta{}_1,\delta{}_2,s-r\}}$ is suitable for $[\,a,t\,]$, because
$x,y\in[\,a,t\,]$ and $|x-y|\loe s-r$ $\Rightarrow$ $x,y\in[\,a,s\,]$ or $x,y\in[\,r,t\,]$.
\smallskip
\CB Let $s\in[\,a,b\,]$. {\sl Case 1:} $s\notin C$. There exists $\delta>0$ with $[\,
s-\delta,s+\delta\,]\cap C=\emptyset$. Then ${[\,a,b\,]}\cap{[\,s-\delta,s+\delta\,]}
\in\cal I$. {\sl Case 2:} $s\in C$. Since $f$ is continuous, there exists $\delta>0$
such that $x\in{[\,s-\delta,s+\delta\,]}\cap C$ $\Rightarrow$ $|f(x)-f(s)|<{\varepsilon
\over2}$. Then ${[\,a,b\,]}\cap{[\,s-\delta,s+\delta\,]}\in\cal I$, because $x,y\in{[\,s-
\delta,s+\delta\,]}\cap C$ $\Rightarrow$ $|f(x)-f(y)|\loe|f(x)-f(s)|+|f(s)-f(y)<\varepsilon$.
\smallskip
By the second key-lemma we conclude that $[\,a,b\,]\in\cal I$, which proves the
theorem.\eop

\bigskip


{\bf Remark.} \ Cousin's theorem was already proved in [5]. For the sake of
completeness, we prove this result once again. At first, we recall a few definitions
coming from the theory of the generalized Riemann integral (see [1] for details).

\bigskip


{\bf Definition.} \ {\bf1)} \ A {\sl gauge} on the interval $[\,a,b\,]$ is a positive
function $\delta:[\,a,b\,]\to\RN$.
\smallskip
{\bf2)} \ A {\sl tagged partition} of $[\,a,b\,]$ is a finite sequence $a=a_0<a_1<
\dots<a_n=b$ together with numbers $x_i\in[\,a_{i-1},a_i\,]$ for all $i=1,\dots,n$.
\smallskip
{\bf3)} \ Given a gauge $\delta:[\,a,b\,]\to\RN$, we say that a tagged partition is
{\sl$\delta$-fine} if it satisfies the condition $[\,a_{i-1},a_i\,]\subset[\,x_i-\delta
(x_i),x_i+\delta(x_i)\,]$ for all $i=1,\dots,n$.
\smallskip
$\delta$-fine tagged partitions will be called $\delta$-fine partitions for short.

\bigskip


\line{{\bf Cousin's theorem.} \ For every gauge $\delta$ on $[\,a,b\,]$, there exists
a $\delta$-fine partition of $[\,a,b\,]$.}

\medskip


{\bf Proof.} \ We consider the set ${\cal I}=\{I\subset[\,a,b\,]\mathrel/\hbox{there
exists a $\delta$-fine partition of }I\}$.
\smallskip
\CA Let $[\,a,s\,]\in\cal I$ and $[\,s,t\,]\in\cal I$. There exists two $\delta$-fine
partitions ${\cal P}_1$ of $[\,a,s\,]$ and ${\cal P}_2$\break of $[\,s,t\,]$. Then
${\cal P}_1\cup{\cal P}_2$ is a $\delta$-fine partition of $[\,a,t\,]$.
\smallskip
\CL Let $s\in{]\,a,b\,]}$ . We put $r=\max{\{a,s-\delta(s)\}}$. Then for every $r\loe
x<s$, the interval $[\,x,s\,]$ together with the tag $s$ is a $\delta$-fine partition
of $[\,x,s\,]$, and hence $[\,x,s\,]\in\cal I$.
\smallskip
The proof of (r) is similar. By the first key-lemma we conclude that $[\,a,b\,]\in
\cal I$.\eop

\bigskip


{\bf Remark.} \ The idea of the following proofs can be found in [12].

\bigskip


{\bf Dini's theorem.} \ Let $(f_n)$ be a decreasing sequence of upper
semicontinuous functions defined on a closed subset $C\subset[\,a,b\,]$.
We suppose that this sequence converges pointwise to zero. Then the
convergence is uniform.

\medskip


{\bf Proof.} \ Let $\varepsilon>0$. We consider the set
\medskip
\centerline{${\cal I}=\{I\subset[\,a,b\,]\mathrel/\hbox{there exists }n\in\NN$ such
that $0\loe f_n(x)<\varepsilon$ for all $x\in I\cap C\}$}
\medskip
\CO Let $[\,a,s\,]\in\cal I$ and $[\,r,t\,]\in\cal I$ with $r<s\loe t$. There exist a
suitable $m$ for $[\,a,s\,]$ and\break a suitable $n$ for $[\,r,t\,]$. Then $p=\max
\{m,n\}$ is suitable for $[\,a,t\,]$.
\smallskip
\CB Let $s\in[\,a,b\,]$. {\sl Case 1:} $s\notin C$. There exists $\delta>0$ with
$[\,s-\delta,s+\delta\,]\cap C=\emptyset$. Then ${[\,a,b\,]}\cap{[\,s-\delta,s+\delta
\,]}\in\cal I$. {\sl Case 2:} $s\in C$. There exists $n\in\NN$ with $0\loe f_n(s)<
\varepsilon$. And by semicontinuity,  there exists $\delta>0$ such that $0\loe
f_n(x)<\varepsilon$ for all $x\in[\,s-\delta,s+\delta\,]\cap C$. Then $\abc\in\cal I$.
\smallskip
By the second key-lemma we conclude that $[\,a,b\,]\in\cal I$, which proves the
theorem.\eop

\bigskip


{\bf Heine-Borel theorem.} \ Let $C\subset[\,a,b\,]$ be a closed subset. Then
any open cover $\cal O$ of the subset $C$ has a finite subcover.

\medskip


{\bf Proof.} \ We consider the set ${\cal I}=\{I\subset[\,a,b\,]\mathrel/I\cap C$
has a finite subcover$\}$.
\smallskip
\CO Let $[\,a,s\,]\in\cal I$ and $[\,r,t\,]\in\cal I$ with $r<s\loe t$. There exist two
finite subcovers ${\cal O}_1$ of $[\,a,s\,]\cap C$ and ${\cal O}_2$ of $[\,r,t\,]\cap
C$. Then ${\cal O}_1\cup{\cal O}_2$ is a finite subcover of $[\,a,t\,]\cap C$.
\smallskip
\CB Let $s\in[\,a,b\,]$. {\sl Case 1:} $s\notin C$. There exists $\delta>0$ with
$[\,s-\delta,s+\delta\,]\cap C=\emptyset$. Then $\emptyset$ is a finite subcover,
and $\abc\in\cal I$. {\sl Case 2:} $s\in C$. Since $\cal O$ is a cover, we have
$s\in O$ for some $O\in\cal O$. There exists $\delta>0$ with $[\,s-\delta,s+
\delta\,]\subset O$. Then $\{O\}$ is\break a finite subcover, and $\abc\in\cal I$.
\smallskip
By the second key-lemma we conclude that $[\,a,b\,]\in\cal I$.\eop

\bigskip


{\bf Bolzano-Weierstrass theorem.} \ If $F\subset[\,a,b\,]$ is a subset with no
accumulation point, then $F$ is finite.

\smallskip


{\bf Proof.} \ We consider the set ${\cal I}=\{I\subset[\,a,b\,]\mathrel/I\cap F$ is
finite$\}$.
\smallskip
\CO If $[\,a,s\,]\cap F$ and $[\,r,t\,]\cap F$ are finite, then $\bigl([\,a,s\,]\cup[\,r,
t\,]\bigr)\cap F$ is also finite.
\smallskip
\CB Let $s\in[\,a,b\,]$. By hypothesis, $s$ is not an accumulation point of the
subset $F$. Then there exists $\delta>0$ with $[\,s-\delta,s+\delta\,]\cap F
\subset\{s\}$, and hence $\abc\in\cal I$.
\smallskip
By the second key-lemma we conclude that $[\,a,b\,]\in\cal I$.\eop

\bigskip


{\bf Cantor's intersection theorem.} \ Let $(C_n)$ be a decreasing sequence of
closed subsets of $[\,a,b\,]$ whose intersection is empty. Then there exists $n
\in\NN$ with $C_n=\emptyset$.

\medskip


{\bf Proof.} \ We consider the set ${\cal I}=\{I\subset[\,a,b\,]\mathrel/\hbox{there
exists }n\in\NN$ with $I\cap C_n=\emptyset\}$.
\smallskip
\CO Let $[\,a,s\,]\in\cal I$ and $[\,r,t\,]\in\cal I$ with $r<s\loe t$. There exist a
suitable $m$ for $[\,a,s\,]$ and\break a suitable $n$ for $[\,r,t\,]$. Then $p=\max
{\{m,n\}}$ is suitable for $[\,a,t\,]$.
\smallskip
\CB Let $s\in[\,a,b\,]$. According to the hypothesis, we have $s\notin C_n$ for
some $n\in\NN$. There exists $\delta>0$ with $[\,s-\delta,s+\delta\,]\cap C_n=
\emptyset$. Then $\abc\in\cal I$.
\smallskip
By the second key-lemma we conclude that $[\,a,b\,]\in\cal I$.\eop

\bigskip


{\bf Historical note.} \ The idea of unifying proofs using some interval induction is
very old. It already appeared in the 1950s in [6], [4], [8], [12] and [14]. The main
lemma in [14] is close to the second key-lemma. Thomson's Lemma [15] on full
covering, which is another version of Cousin's theorem, is used in [2], [3], [9] and
[16]. Finally, tagged partitions are introduced in [7] and [13]. But this was already
the case twenty years earlier in the book [10] (in French). An interesting
bibliography can be found in [11].

\bigskip


\setbox7=\hbox{[0] \ \ }

{\bf References}

\smallskip

$\hbox to\wd7{[1]\hfill}$R.\ G.\ Bartle and D.\ R.\ Sherbert, {\sl Introduction to real
analysis 4th edition,}\par
${}\kern\wd7$New York, John Wiley \& Sons, 2011

\smallskip

$\hbox to\wd7{[2]\hfill}$M.\ W.\ Botsko, {\sl A unified treatment of various theorems in
elementary analysis}
\par
${}\kern\wd7$Amer.\ Math.\ Monthly {\bf94} (1987), 450--452

\smallskip

$\hbox to\wd7{[3]\hfill}$M.\ W.\ Botsko, {\sl The use of full covers in real analysis}
\par
${}\kern\wd7$Amer.\ Math.\ Monthly {\bf96} (1989), 328--333

\smallskip

$\hbox to\wd7{[4]\hfill}$W.\ L.\ Duren, {\sl Mathematical induction in sets}
\par
${}\kern\wd7$Amer.\ Math.\ Monthly {\bf64} (1957), 19--22

\smallskip

$\hbox to\wd7{[5]\hfill}$C.-A.\ Faure, {\sl The last proof of extreme value theorem
and intermediate value theorem}
\par
${}\kern\wd7$arXiv:2209.12682 [math.HO]

\smallskip

$\hbox to\wd7{[6]\hfill}$L.\ R.\ Ford, {\sl Interval additive propositions}
\par
${}\kern\wd7$Amer.\ Math.\ Monthly {\bf64} (1957), 106--108

\smallskip

$\hbox to\wd7{[7]\hfill}$R.\ A.\ Gordon, {\sl The use of tagged partitions in elementary
real analysis}
\par
${}\kern\wd7$Amer.\ Math.\ Monthly {\bf105} (1998), 107--117 \& 886

\smallskip

$\hbox to\wd7{[8]\hfill}$G.\ Jungck, {\sl Interval induction}
\par
${}\kern\wd7$Amer.\ Math.\ Monthly {\bf73} (1966), 295--297

\smallskip

$\hbox to\wd7{[9]\hfill}$K.\ Klaimon, {\sl More applications of full covering}
\par
${}\kern\wd7$Pi Mu Epsilon J.\ {\bf9} (1990), 156-161

\smallskip

$\hbox to\wd7{[10]\hfill}$J.\ Mawhin, {\sl Introduction \`a l'analyse 1\`ere \'edition}
\par
${}\kern\wd7$Louvain-la-Neuve, Cabay, 1979

\smallskip

$\hbox to\wd7{[11]\hfill}$J.\ Mawhin, {\sl Initiation \`a la compacit\'e : une variante} 
\par
\line{${}\kern\wd7$L'enseignement de l'analyse aux d\'ebutants, Louvain-la-Neuve,
Erasme, 1992, 109--125}

\smallskip

$\hbox to\wd7{[12]\hfill}$R.\ M.\ F.\ Moss, G.\ T.\ Roberts, {\sl A creeping lemma}
\par
${}\kern\wd7$Amer.\ Math.\ Monthly {\bf75} (1968), 649--652

\smallskip

$\hbox to\wd7{[13]\hfill}$S.\ Prongjit and W.\ Sodsiri, {\sl Applications of $\delta$-fine
tagged partitions in real analysis}
\par
${}\kern\wd7$Far East J.\ Math.\ Sci.\ {\bf91} (2014), 97--109 \& 233--245

\smallskip

$\hbox to\wd7{[14]\hfill}$P.\ Shanahan, {\sl A unified proof of several basic theorems
of real analysis}
\par
${}\kern\wd7$Amer.\ Math.\ Monthly {\bf79} (1972), 895--898 \& {\bf81} (1974), 890--891

\smallskip

$\hbox to\wd7{[15]\hfill}$B.\ S.\ Thomson, {\sl On full covering properties}
\par
${}\kern\wd7$Real Anal.\ Exchange {\bf6} (1980--81), 77--93

\smallskip

$\hbox to\wd7{[16]\hfill}$K.\ Zangara and J.\ Marafino, {\sl Applications of full covers
in real analysis}
\par
${}\kern\wd7$Involve {\bf2} (2009), 297--304

\bigskip\medskip


\centerline{\vbox{\hbox{Gymnase de la Cit\'e, place de la Cath\'edrale 1, 1014
Lausanne, Switzerland}
\hbox{Email address: {\tt claudealain.faure@eduvaud.ch}}}}

\end